\documentclass[11pt,twoside,reqno,centertags,draft]{amsart}
\usepackage{amsfonts}
\usepackage{color,enumitem,graphicx}
\usepackage[colorlinks=true,urlcolor=blue,
citecolor=red,linkcolor=blue,linktocpage,pdfpagelabels,
bookmarksnumbered,bookmarksopen]{hyperref}

  \usepackage{amsmath,amsthm,amsfonts,amssymb}

\begin{document}
\title{ Existence and nonexistence of solutions to \\
Choquard equations}
\date{}
\maketitle

\vspace{ -1\baselineskip}

{\small
\begin{center}

\medskip

  {\sc  Wanwan Wang}
  \medskip

Department of Mathematics, Jiangxi Normal University,\\
 Nanchang, Jiangxi 330022, PR China
\end{center}
}

\renewcommand{\thefootnote}{}

\footnote{E-mail address: wwwang2014@yeah.net (W. Wang).}
\footnote{MSC2010: 35B09, 35B33,  35J61.}
\footnote{Keywords: Choquard equation;  Riesz potential; Poho\v{z}aev identity;   Existence;  Nonexistence.}

\begin{quote}
{\bf Abstract.} In this paper, we establish the existence of ground state solutions for Choquard equations
\begin{equation}\label{eq 1}
 - \Delta u + u = q\,(I_\alpha \ast |u|^p) |u|^{q - 2} u+p\,(I_\alpha \ast |u|^q) |u|^{p - 2} u\quad {\rm in }\quad  \mathbb{R}^N,
\end{equation}
where $N \ge 3$,   $\alpha \in (0, N)$,
  $I_\alpha: \mathbb{R}^N \to \mathbb{R}$ is the Riesz potential,  $p,\,q >0$ satisfying that
\begin{equation}\label{eq 2}
\frac{2(N+\alpha)}{N}<p+q< \frac{2(N+\alpha)}{N-2}.
\end{equation}
Moreover,  we prove a Poho\v{z}aev type identity for problem (\ref{eq 1}), which implies the non-existence result for the problem when $(p,q)$ does not satisfy the condition  (\ref{eq 2}).

\end{quote}

\medskip

\newcommand{\N}{\mathbb{N}}
\newcommand{\R}{\mathbb{R}}
\newcommand{\Z}{\mathbb{Z}}

\newcommand{\cA}{{\mathcal A}}
\newcommand{\cB}{{\mathcal B}}
\newcommand{\cC}{{\mathcal C}}
\newcommand{\cD}{{\mathcal D}}
\newcommand{\cE}{{\mathcal E}}
\newcommand{\cF}{{\mathcal F}}
\newcommand{\cG}{{\mathcal G}}
\newcommand{\cH}{{\mathcal H}}
\newcommand{\cI}{{\mathcal I}}
\newcommand{\cJ}{{\mathcal J}}
\newcommand{\cK}{{\mathcal K}}
\newcommand{\cL}{{\mathcal L}}
\newcommand{\cM}{{\mathcal M}}
\newcommand{\cN}{{\mathcal N}}
\newcommand{\cO}{{\mathcal O}}
\newcommand{\cP}{{\mathcal P}}
\newcommand{\cQ}{{\mathcal Q}}
\newcommand{\cR}{{\mathcal R}}
\newcommand{\cS}{{\mathcal S}}
\newcommand{\cT}{{\mathcal T}}
\newcommand{\cU}{{\mathcal U}}
\newcommand{\cV}{{\mathcal V}}
\newcommand{\cW}{{\mathcal W}}
\newcommand{\cX}{{\mathcal X}}
\newcommand{\cY}{{\mathcal Y}}
\newcommand{\cZ}{{\mathcal Z}}

\newcommand{\abs}[1]{\lvert#1\rvert}
\newcommand{\xabs}[1]{\left\lvert#1\right\rvert}
\newcommand{\norm}[1]{\lVert#1\rVert}

\newcommand{\loc}{\mathrm{loc}}
\newcommand{\p}{\partial}
\newcommand{\h}{\hskip 5mm}
\newcommand{\ti}{\widetilde}
\newcommand{\D}{\Delta}
\newcommand{\e}{\epsilon}
\newcommand{\bs}{\backslash}
\newcommand{\ep}{\emptyset}
\newcommand{\su}{\subset}
\newcommand{\ds}{\displaystyle}
\newcommand{\ld}{\lambda}
\newcommand{\vp}{\varphi}
\newcommand{\wpp}{W_0^{1,\ p}(\Omega)}
\newcommand{\ino}{\int_\Omega}
\newcommand{\bo}{\overline{\Omega}}
\newcommand{\ccc}{\cC_0^1(\bo)}
\newcommand{\iii}{\opint_{D_1}D_i}

\numberwithin{equation}{section}

\vskip 0.2cm \arraycolsep1.5pt
\newtheorem{lemma}{Lemma}[section]
\newtheorem{theorem}{Theorem}[section]
\newtheorem{definition}{Definition}[section]
\newtheorem{proposition}{Proposition}[section]
\newtheorem{remark}{Remark}[section]
\newtheorem{corollary}{Corollary}[section]

\setcounter{equation}{0}
\section{Introduction}

This paper is devoted to the study of existence results for nonnegative solutions  of Choquard equations
\begin{equation}\label{eq 1.1}
 - \Delta u + u = q(I_\alpha \ast |u|^p) |u|^{q - 2} u+p(I_\alpha \ast |u|^q) |u|^{p - 2} u\quad {\rm in }\ \   \R^N,\quad  \  u\in H^1(\R^N),
\end{equation}
 where  $p,\,q>0$, $N \ge 3$,   $\alpha \in (0, N)$
 and
  $I_\alpha: \R^N \to \R$ is the Riesz potential defined by
$$
  I_\alpha(x)= \frac{\Gamma(\frac{N-\alpha}{2})|x|^{\alpha-N}}
                   {\pi^{N/2}2^{\alpha}\Gamma(\frac{\alpha}{2})}
$$
with $\Gamma$ being the Gamma function, see \cite{RM}.

As early as in 1954, the Choquard equation
$$
\left\{ \arraycolsep=1pt
\begin{array}{lll}
 -\Delta u+u=(I_2\ast|u|^2)u\ \ \ &
{\rm in} \quad \R^3,\\[2mm]
\lim_{|x|\to+\infty}u(x)=0
\end{array}
\right.
$$
has appeared in the context of various physical models. It seems to originate
from H. Fr\"{o}hlich and S. Pekar＊s model of the polaron, where free electrons
in an ionic lattice interact with phonons associated to deformations of the
lattice or with the polarisation that it creates on the medium (interaction
of an electron with its own hole) \cite{F,F1,P}. The Choquard equation was
also introduced by Ph. Choquard in 1976 in the modelling of a one-component
plasma.

The existence and qualitative properties of solutions of Choquard  equations
 have been widely studied in the last decades. In \cite{L1}, Lieb proved the
 existence and uniqueness, up to translations, of the ground state. Later on,
 in \cite{L2}, Lions showed the existence of a sequence of radially symmetric solutions.
 In \cite{CCS,G1,G2,Le,MZ} the authors considered the regularity, positivity and radial symmetry of
 the ground states and derived decay property at infinity as well.
 Especially, Moroz and Van Schaftingen in \cite{VJ} studied the  generalized  Choquard equation
\begin{equation}\label{eq 10}
\left\{ \arraycolsep=1pt
\begin{array}{lll}
 -\Delta u+u=(I_\alpha\ast|u|^p)|u|^{p-2}u\ \ \ &
{\rm in} \quad \R^N,\\[2mm]
\lim_{|x|\to+\infty}u(x)=0,
\end{array}
\right.
\end{equation}
they showed that solutions of problem (\ref{eq 10}) are, at least formally, critical points of the  functional
 $F$ defined for a function $u:\R^N\to\R$ by
 $$F(u)=\frac12\int_{\R^N}(|\nabla u(x)|^2+|u(x)|^2)\,dx-\frac1{2p} \int_{\mathbb{R}^N}(I_\alpha*|u|^p)(x)|u(x)|^p\,dx.$$

In the present paper, we are interested in studying the existence of ground states solutions for Choquard problem  (\ref{eq 1.1}).
We note that problem  (\ref{eq 1.1}) has a variational structure: the critical points of the function
 $$E(u)=\frac12\int_{\R^N}(|\nabla u(x)|^2+|u(x)|^2)\,dx- \int_{\mathbb{R}^N}(I_\alpha*|u|^p)(x)|u(x)|^q\,dx$$
 are solutions of (\ref{eq 1.1}). We state the existence results as follows.

\begin{theorem} \label{theorem 1}
Suppose that $N \ge 3,\,  \alpha \in (0, N)$ and $p,\,q >0$ satisfying that
\begin{equation}\label{e}
\frac{2(N+\alpha)}{N}<p+q< \frac{2(N+\alpha)}{N-2}.
\end{equation}
Then problem (\ref{eq 1.1}) admits  a positive ground state solution.
\end{theorem}

To prove the existence result in Theorem \ref{theorem 1}, we apply the critical points theory to the associated minimizing problem
\begin{equation}\label{4e}
M_p=\inf\left\{\int_{\mathbb R^N}(|\nabla u(x)|^2+|u(x)|^2)\,dx: \int_{\mathbb{R}^N}(I_\alpha*|u|^p)(x)|u(x)|^q\,dx=1\right\}.
\end{equation}
By Hardy-Littlewood-Sobolev inequality, which states that
if $t,r>1$ and $\frac{1}{t}+\frac{N-\alpha}{N}+\frac{1}{r}=2$,
$f\in L^t(\R^N)$ and $h\in L^r(\R^N)$, then there exists a sharp constant $C(t,N,\alpha,r)$, independent of $f,h$, such that
$$\int_{\R^N}\int_{\R^N}\frac{f(x)h(y)}{|x-y|^{N-\alpha}}\, dxdy \leq C(t,N,\alpha,r)\|f\|_{L^t(\R^N)}\|h\|_{L^r(\R^N)},$$
see for instance \cite{G2,VJ}, we see that  $M_p>0$.
Then we use the  concentration compactness argument and a nonlocal version of Brezis-Lieb lemma to prove that $M_p$ can be achieved.
The minimization of $M_p$ is a  nontrivial solution of
(\ref{eq 1.1}).

The second aim of this paper is to establish the Poho\v{z}aev type identity for  (\ref{eq 1.1}) and  obtain the non-existence results as follows.
\begin{theorem}\label{theorem 2}
Let $u$ be a nonnegative solution of  (\ref{eq 1.1})
with $p,\,q>0$ satisfying that
\begin{equation}\label{n}
{p+q}\ge \frac{2(N+\alpha)}{N-2}\quad  {\rm or} \quad {p+q} \leq \frac{2(N+\alpha)}{N}.
\end{equation}
Assume that $u\in H^1(\mathbb{R}^N)\cap L^\frac{N(p+q)}{N+\alpha}(\mathbb{R}^N)$ and  $|\nabla u|\in H_{loc}^1(\mathbb {R}^N)$.
Then $u$ is a trivial solution of  (\ref{eq 1.1}).
\end{theorem}

The content of the paper is the following: in Section 2 we provide some technical preliminaries;
in Section 3 we prove the existence of ground state solutions of  (\ref{eq 1.1})
 in Theorem \ref{theorem 1} by the critical points theory; in Section 4 we show the Poho\v{z}aev type identity and then
 prove the non-existence results in Theorem \ref{theorem 2}.

\setcounter{equation}{0}
\section{Preliminaries}
The purpose of this section is to introduce some preliminaries.

\begin{lemma}\label{lemma 1}\cite {W}
Let $\Omega $ be a domain in $\mathbb{R}^N$, $t>1$ and $\{w_m\}_{m\in \mathbb{N}}$ be a bounded sequence in $L^s(\Omega)$.
If $w_m\to w$ almost everywhere on $\Omega$ as $m\to \infty$, then for every $r\in[1,s]$, we have that
$$\lim_{m\to\infty}{\int_\Omega|{|w_m|^r-|w_m-w|^r-|w|^r}|^\frac{t}{r}}\,dx=0.$$
 \end{lemma}

\begin{lemma}\label{lemma 2}
Let  $\alpha \in(0,N)$, $\frac{2(N+\alpha)}{N}<p+q< \frac{2(N+\alpha)}{N-2}$ and $\{w_m\}_{m\in \mathbb{N}}$ be a bounded sequence in
$L^{\frac{N(p+q)}{N+\alpha}}(\mathbb{R}^N)$. Assume that

\quad (i) ${w_m}$ weakly converges to ${w}$ in $L^\frac{N(p+q)}{N+\alpha}(\mathbb{R}^N)$;

\quad (ii)  ${w_m}\to w$ almost everywhere on  $\mathbb{R}^N$.\\
Then
\begin{eqnarray*}
&&\lim_{m\to \infty}\left[\int_{\mathbb{R}^N}(I_\alpha*|w_m|^p)(x)|w_m(x)|^q\,dx-\int _{\mathbb{R}^N}(I_\alpha*|w_m-w|^p)(x)|(w_m-w)(x)|^q\,dx\right]
\\&&=\int_{\mathbb{R}^N}(I_\alpha*|w|^p)(x)|w(x)|^q\,dx.
\end{eqnarray*}
 \end{lemma}

\noindent{\bf Proof.} By direct computation,
we have that
 \begin{eqnarray*}
 && \int_{\mathbb{R}^N}(I_\alpha*|w_m|^p)(x)|w_m(x)|^q\,dx-\int _{\mathbb{R}^N}(I_\alpha*|w_m-w|^p)(x)|(w_m-w)(x)|^q\,dx
\\&=&\int_{\mathbb{R}^N}(I_\alpha*(|w_m|^p-|w_m-w|^p))(x)(|w_m(x)|^q-|(w_m-w)(x)|^q)\,dx
\\&&+\int_{\mathbb{R}^N}(I_\alpha*(|w_m|^p-|w_m-w|^p))(x)|(w_m-w)(x)|^q\,dx
\\&&+\int_{\mathbb{R}^N}(I_\alpha*|w_m-w|^p)(x)(|w_m(x)|^q-|(w_m-w)(x)|^q)\,dx
\\&:=&A_1+A_2+A_3.
\end{eqnarray*}
We look at each of these integrals separately. First, we use the H\"{o}lder inequality to obtain that
\begin{eqnarray*}
  A_2&=&\int_{\mathbb{R}^N}(I_\alpha*(|w_m|^p-|w_m-w|^p-|w|^p))(x)|(w_m-w)(x)|^q\,dx
\\&& \qquad \ \  +\int_{\mathbb{R}^N}(I_\alpha*|w|^p)(x)|(w_m-w)(x)|^q\,dx
\\&\leq &\left({\int_{\mathbb{R}^N}|\left(I_\alpha*(|w_m|^p-|w_m-w|^p-|w|^p)\right)|^\frac{N(p+q)}{Np-\alpha q}}(x)\,dx\right)^\frac{Np-\alpha q}{N(p+q)}
\\&&\cdot\left(\int_{\mathbb{R}^N}(|(w_m-w)(x)|^q)^\frac{N(p+q)}{(N+\alpha)q}\,dx\right)^\frac{(N+\alpha)q}{N(p+q)}+\int_{\mathbb{R}^N}(I_\alpha*|w|^p)(x)|(w_m-w)(x)|^q\,dx.
\end{eqnarray*}
Using Lemma \ref{lemma 1} with $r=p$ and $t=\frac{N(p+q)}{N+\alpha}$, we know that $|w_m|^p-|w_m-w|^p\to |w|^p$,
strongly in $L^\frac{N(p+q)}{(N+\alpha)p}(\mathbb{R}^N)$ as $m\to\infty$. By the Hardy-Littlewood-Sobolev inequality, this implies that
$I_\alpha*(|w_m|^p-|w_m-w|^p)\to I_\alpha*|w|^p$ in $L^\frac{N(p+q)}{Np-\alpha q}(\mathbb{R}^N)$ as $m\to\infty$. Since    $|w_m-w|^q\rightharpoonup 0$
in $L^\frac{N(p+q)}{(N+\alpha)q}(\mathbb{R}^N)$ as $m\to\infty$, then $A_2\to0$
as $m\to\infty$. We next deal with the term $A_3$. We observe that
\begin{eqnarray*}
  A_3&=&\int_{\mathbb{R}^N}(I_\alpha*|w_m-w|^p)(x)(|w_m(x)|^q-|(w_m-w)(x)|^q-|w(x)|^q)\,dx
\\&&  \ \  +\int_{\mathbb{R}^N}(I_\alpha*|w_m-w|^p)(x)|w(x)|^q\,dx
\\&\leq &\left({\int_{\mathbb{R}^N}|\left(I_\alpha*(|w_m-w|^p\right)|^\frac{N(p+q)}{Np-\alpha q}}(x)\,dx\right)^\frac{Np-\alpha q}{N(p+q)}
\\&&\cdot\left(\int_{\mathbb{R}^N}(|w_m(x)|^q-|(w_m-w)(x)|^q-|w(x)|^q)^\frac{N(p+q)}{(N+\alpha)q}\,dx\right)^\frac{(N+\alpha)q}{N(p+q)}
\\&& \qquad \qquad \qquad \qquad\qquad \qquad \qquad +\int_{\mathbb{R}^N}(I_\alpha*|w_m-w|^p)(x)|w(x)|^q\,dx,
\end{eqnarray*}
which implies $A_3\to 0$ as  $m\to\infty$ by Lemma \ref{lemma 1}.
Finally, we note that
 $$A_1\to \int_{\mathbb{R}^N}(I_\alpha*|w|^p)(x)|w(x)|^q\,dx$$
 as $m\to\infty$.
The proof ends. \hfill$\Box$

\setcounter{equation}{0}
\section{Ground state solutions}

In this section, we establish the existence of ground state solutions of  (\ref{eq 1.1}). Let us consider the  minimizing problem
\begin{equation}\label{4e}
M_p=\inf\left\{\int_{\mathbb R^N}(|\nabla u(x)|^2+|u(x)|^2)\,dx: \int_{\mathbb{R}^N}(I_\alpha*|u|^p)(x)|u(x)|^q\,dx=1\right\},
\end{equation}
defined on ${H^1(\mathbb{R}^N)}$. By Hardy-Littlewood-Sobolev inequality, we note that  $M_p$ is well defined.
\begin{proposition}\label{proposition 3.1}
The minimizing problem $M_p$ is achieved by a function $v\in {H^1(\mathbb{R}^N)}$, which is a solution of (\ref{eq 1.1})
up to a translation.
\end{proposition}
We will use the concentration-compactness principle \cite{L} to prove Proposition \ref{proposition 3.1}.
To this end, we introduce the following vanishing type lemma. Let  $B_r(x)$ denote the ball centered at $x\in \mathbb{R}^N$ with radius r.
\begin{lemma}\label{lemma 3}
Let $2\leq s<2^*=\frac{2N}{N-2}$ and $r>0$. Suppose that $\{v_m\}_{m\in\mathbb{N}}$ is a bounded sequence in $H^1(\mathbb{R}^N)$
and
$$\sup_{z\in{\mathbb{R}^N}}\int _{B_r(z)}|v_m(x)|^s\,dx\to 0$$
as $m\to\infty$.
Then  for $\frac{(N+\alpha)s}{N}<p+q<\frac{2(N+\alpha)}{N-2}$, we have that
 $$\int_{\mathbb{R}^N}(I_\alpha*|v_m|^p)(x)|v_m(x)|^q\,dx\to0$$
as $m\to\infty$.
\end{lemma}

\noindent {\bf Proof.}
Let  $l=\frac{p+q}{q}\frac{N}{N+\alpha}$ and $t=\frac{p+q}{p}\frac{N}{N+\alpha}$, then $lq=pt$,  by Hardy-Littlewood-Sobolev inequality, there exists $C>0$ such that
\begin{eqnarray*}
 && \int_{\mathbb{R}^N}(I_\alpha*|v_m|^p)(x)|v_m(x)|^q\,dx\ =\int_{\mathbb{R}^N}\int_{\mathbb{R}^N}\frac{|v_m(x)|^q|v_m(z)|^p}{|x-z|^{N-\alpha}}\,dxdz
\\& & \leq  C\|{|v_m|^q}\|_{L^l(\mathbb{R}^N)}\|{|v_m|^p}\|_{L^t(\mathbb{R}^N)} \  =C\left(\int_{\mathbb{R}^N}|v_m(x)|^\frac{N(p+q)}{N+\alpha}\,dx\right)^\frac{N+\alpha}{N}.
\end{eqnarray*}
Since $s<\frac{N(p+q)}{N+\alpha}<2^*$, using the classical Vanishing  Lemma (see Lemma 1.21 in \cite{W}), it is true that $v_m\to0$ in
$L^{\frac{N(p+q)}{N+\alpha}}(\mathbb{R}^N)$ as
$m\to\infty$.
Thus,  $$\int_{\mathbb{R}^N}(I_\alpha*|v_m|^p)(x)|v_m(x)|^q\,dx\to0$$
as $m\to\infty$.
The proof is complete.\hfill$\Box$

We now  prove proposition \ref{proposition 3.1}.

\noindent{\bf Proof of Proposition \ref{proposition 3.1}}
\ Let $\{v_m\}_{m\in \mathbb{N}}\subset H^1(\mathbb{R}^N)$ be a minimizing sequence of $M_p$ and satisfy that
$$\int_{\mathbb{R}^N}(I_\alpha*|v_m|^p)(x)|v_m(x)|^q\,dx=1$$
and
$$\int_{\mathbb{R}^N}(|\nabla v_m(x)|^2+|v_m(x)|^2)\,dx \to M_p$$
as $m\to\infty$.

By Lemma \ref{lemma 3}, there exists $\delta>0$ such that
$$\delta=\liminf_{m\to\infty}\sup_{z\in\mathbb{R}^N}\int_{B_1(z)}|v_m(x)|^2\,dx>0.$$
Going if necessary to a subsequence, we may assume the existence of  $\{z_m\}_{m\in\mathbb{N}}\in\mathbb{R}^N$ such that
$$\int_{B_1(z_m)}|v_m(x)|^2\,dx>\frac{\delta}{2}.$$
Let us denote $w_m(x)=v_m(x-z_m)$, then we have that
$$\int_{\mathbb{R}^N}(I_\alpha*|w_m|^p)(x)|w_m(x)|^q\,dx=1,\quad \ \ \int_{\mathbb R^N}(|\nabla w_m(x)|^2+|w_m(x)|^2)\,dx \to M_p$$
and
\begin{equation}\label{eq 32}
\int_{B_1(0)}|w_m(x)|^2\,dx>\frac{\delta}{2}.
\end{equation}
Since $\{w_m\}_m\in\mathbb{N}$ is bounded in $H^1(\mathbb{R}^N)$,   there exists $w$ such that
  $w_m\rightharpoonup w$ in $H^1(\R^N)$, $w_m\to w$ in $L_{loc}^2(\R^N)$ and $w_m \to w$ almost everywhere on $\mathbb{R}^N$.
  Combining with (\ref{eq 32}), we have that  $w\neq 0$ almost everywhere on $\mathbb{R}^N$.
Then $\int_{\mathbb{R}^N}(I_\alpha*|w|^p)(x)|w(x)|^q\,dx\neq 0$.

Using  Lemma \ref{lemma 2}, we obtain that
\begin{eqnarray*}
1= \int_{\mathbb{R}^N}(I_\alpha*|w|^p)(x)|w(x)|^q\,dx+\lim_{m\to \infty}\int_{\mathbb{R}^N}(I_\alpha*|w_m-w|^p)(x)|(w_m-w)(x)|^q\,dx
\end{eqnarray*}
and
\begin{eqnarray*}
&&M_p=\lim_{m\to\infty}{\|w_m\|_{H^1(\mathbb{R}^N)}^2}={\|w\|_{H^1(\mathbb{R}^N)}^2}+\lim_{m\to\infty}{\|w_m-w\|_{H^1(\mathbb{R}^N)}^2}
\\&&\geq M_p\left(\int_{\mathbb{R}^N}(I_\alpha*|w|^p)(x)|w(x)|^q\,dx\right)^\frac{2}{p+q}
\\&& \ \quad + M_p\left(\lim_{m\to\infty}\int_{\mathbb{R}^N}(I_\alpha*|w_m-w|^p)(x)|(w_m-w)(x)|^q\,dx\right)^\frac{2}{p+q}
\\&&= M_p \left(\int_{\mathbb{R}^N}(I_\alpha*|w|^p)(x)|w(x)|^q\,dx\right)^\frac{2}{p+q}+M_p\left(1-\int_{\mathbb{R}^N}(I_\alpha*|w|^p)(x)|w(x)|^q\,dx\right)^\frac{2}{p+q}
\end{eqnarray*}
Then  $\int_{\mathbb{R}^N}(I_\alpha*|w|^p)(x)|w(x)|^q\,dx=1$. As a  consequent, we get that $M_p=\|w\|_{H^1(\mathbb {R}^N)}^2$.\\
The proof is completed. \hfill$\Box$

\setcounter{equation}{0}
\section{Nonexistence}

In this section, we prove a Poho\v{z}aev type identity for  (\ref{eq 1.1}), then we obtain the non-existence result of  (\ref{eq 1.1}) by this Poho\v{z}aev type identity.
\begin{lemma}\label{lemma 4}
Let $u\in H^1(\mathbb{R}^N)\cap L^\frac{N(p+q)}{N+\alpha}(\mathbb{R}^N)$ be a solution of  (\ref{eq 1.1}) and $|\nabla u|\in H_{loc}^1(\mathbb {R}^N)$.
Then
\begin{equation}\label{eq 123}
\frac{N-2}{2}\int_{\mathbb{R}^N}|\nabla u(x)|^2\,dx+\frac{N}{2}\int_{\mathbb{R}^N}|u(x)|^2\,dx=(N+\alpha)\int_{\mathbb{R}^N}(I_\alpha*|u|^p)(x)|u(x)|^q\,dx.
\end{equation}
\end{lemma}
\noindent {\bf Proof.}
Let   $\lambda\in(0,\infty)$,  $x\in{\mathbb{R}^N}$ and $\varphi\in C_c^1(\mathbb{R}^N)$ such that $\varphi=1$ in $B_1(0)$, we denote
\begin{equation}\label{eq10}
v_\lambda(x)=\varphi(\lambda x)x\cdot \nabla u(x).
\end{equation}
 Using $v_\lambda$ as a test function in the equation (\ref{eq 1.1}), we find that
\begin{eqnarray*}
\int_{\mathbb{R}^N}\nabla u\cdot\nabla v_\lambda\,dx+\int_{\mathbb{R}^N}u\, v_\lambda\,dx =\int_{\mathbb{R}^N}(q(I_\alpha \ast |u|^p) |u|^{q - 2} u\, v_\lambda
+p(I_\alpha \ast |u|^q) |u|^{p - 2} u \, v_\lambda)\,dx.
\end{eqnarray*}
We look at each of these integrals separately.
Since  $|\nabla u|\in H_{loc}^1(\mathbb {R}^N)$, combining with (\ref{eq10}),
 we have that
\begin{eqnarray*}
\int_{\mathbb{R}^N}\nabla u\cdot\nabla v_\lambda\,dx
=-\int_{\mathbb{R}^N}\left((N-2)\varphi(\lambda x)+\lambda x\cdot \nabla\varphi(\lambda x)\right)\frac{|\nabla u(x)|^2}{2}\,dx,
\end{eqnarray*}
 then
$$\lim_{\lambda\to 0}\int_{\mathbb{R}^N}\nabla u  \cdot \nabla v_\lambda\,dx=-\frac{N-2}{2}\int_{\mathbb{R}^N}|\nabla u|^2\,dx.$$
By the definition of $v_\lambda$, we also can get that
\begin{eqnarray*}
\int_{\mathbb{R}^N}u\, v_\lambda\,dx
&&=\int_{\mathbb{R}^N}u(x)\varphi(\lambda x)x\cdot \nabla u(x)\,dx
=\int_{\mathbb{R}^N}\varphi(\lambda x)x\cdot \nabla(\frac{|u(x)|^2}{2})\,dx
\\&&=-\int_{\mathbb{R}^N}\left(N\varphi(\lambda x)+\lambda x\cdot\nabla\varphi(\lambda x)\right)(\frac{|u(x)|^2}{2})\,dx,
\end{eqnarray*}
by Lebesgue's dominated convergence theorem, it holds
$$\lim_{\lambda\to 0}\int_{\mathbb{R}^N}u\, v_\lambda\,dx=-\frac{N}{2}\int_{\mathbb{R}^N}|u|^2\,dx.$$
Finally, by direct compute, we have that
\begin{eqnarray*}
&&\int_{\mathbb{R}^N}[q(I_\alpha \ast |u|^p) |u|^{q - 2} u v_\lambda+p(I_\alpha \ast |u|^q) |u|^{p - 2} w v_\lambda]\,dx
\\&=&\int_{\mathbb{R}^N}\int_{\mathbb{R}^N}\left(I_\alpha(x-y)\varphi(\lambda x)x\right)\left[|u(y)|^p\nabla(|u(x)|^q)+|u(y)|^q\nabla(|u(x)|^p)\right]\,dxdy
\\&=&\int_{\mathbb{R}^N}\int_{\mathbb{R}^N}I_\alpha(x-y)\left(|u(y)|^p\varphi(\lambda x)x\cdot\nabla(|u(x)|^q)+|u(x)|^q\varphi(\lambda y)y\cdot\nabla(|u(y)|^p)\right)\,dxdy
\\&=&-\int_{\mathbb{R}^N}\int_{\mathbb{R}^N}|u(y)|^p|u(x)|^q\left[I_\alpha(x-y)\left(\lambda\nabla\varphi(\lambda x)x+N\varphi(\lambda x)\right)-\frac{(x-y)\cdot x\varphi(\lambda x)(N-\alpha)}{|x-y|^{N-\alpha+2}}\right.
\\&&\ \left. +I_\alpha(x-y)\left(\lambda\nabla\varphi(\lambda y)y+N\varphi(\lambda y)\right)+\frac{(x-y)\cdot y\varphi(\lambda y)(N-\alpha)}{|x-y|^{N-\alpha+2}}\right]\,dxdy
\end{eqnarray*}
and then
\begin{eqnarray*}
&&\lim_{\lambda\to 0}\int_{\mathbb{R}^N}q(I_\alpha \ast |u|^p) |u|^{q - 2} u \, v_\lambda+p(I_\alpha \ast |u|^q) |u|^{p - 2} u \, v_\lambda\,dx
\\&=&-\int_{\mathbb{R}^N}\int_{\mathbb{R}^N}|u(y)|^p|u(x)|^q\left[2N\cdot I_\alpha(x-y)-(N-\alpha)\frac{(x-y)\cdot(x-y)}{|x-y|^{N-\alpha+2}}\right]\,dxdy
\\&=&-\int_{\mathbb{R}^N}\int_{\mathbb{R}^N}2N\cdot\frac{|u(y)|^p|u(x)|^q}{|x-y|^{N-\alpha}}\,dxdy
+\int_{\mathbb{R}^N}\int_{\mathbb{R}^N}(N-\alpha)\cdot\frac{|u(y)|^p|u(x)|^q}{|x-y|^{N-\alpha}}\,dxdy
\\&=&-\int_{\mathbb{R}^N}\int_{\mathbb{R}^N}(N+\alpha)\cdot\frac{|u(y)|^p|u(x)|^q}{|x-y|^{N-\alpha}}\,dxdy
\\&=&-(N+\alpha)\int_{\mathbb{R}^N}\int_{\mathbb{R}^N}(I_\alpha \ast |u|^p) |u|^q\,dx.
\end{eqnarray*}
The proof ends. \hfill$\Box$

We now prove the nonexistence result in Theorem \ref{theorem 2} by Lemma \ref{lemma 4}.\\

\noindent{\bf Proof of Theorem \ref{theorem 2}.}
Since $u$ is  a solution of problem (\ref{eq 1.1}), then
$$\int_{\mathbb{R}^N}|\nabla u|^2\, dx+\int_{\mathbb{R}^N}|u|^2\,dx =(p+q) \int_{\mathbb{R}^N}(I_\alpha \ast |u|^p) |u|^q \,dx,$$
combining with the  Poho\v{z}aev type identity (\ref{eq 123}), we have that
$$(\frac{N-2}{2}-\frac{N+\alpha}{p+q})\int_{\mathbb{R}^N}|\nabla u|^2\,dx+(\frac{N}{2}-\frac{N+\alpha}{p+q})\int_{\mathbb{R}^N}|u|^2\,dx=0.$$
When
$${p+q} \ge\frac{2(N+\alpha)}{N-2}\qquad {\rm or} \qquad {p+q} \le \frac{2(N+\alpha)}{N},$$
 it holds that $u = 0$. \hfill$\Box$

\bigskip

\noindent{\bf Acknowledgements:} The author would like to express the  warmest gratitude to Prof. Jianfu Yang,
for proposing the problem and for its active participation.
This work is supported by the Jiangxi Provincial Natural Science Foundation (20161ACB20007).

\end{document}